\documentclass[11pt]{article}
\usepackage{german}
\usepackage[T1]{fontenc}
\usepackage[ansinew]{inputenc}
\usepackage{latexsym}
\usepackage{textcomp}
\usepackage{color}
\usepackage{amscd}
\usepackage{amsmath}

\input{amssym.def}
\input{amssym}

     \newcommand{\fa}{\goth{a}}
     \newcommand{\fb}{\goth{b}}

     \newcommand{\ff}{\goth{f}}

     \newcommand{\fo}{\goth{o}}
     
     \newcommand{\fq}{\goth{q}}
     
     \newcommand{\fs}{\goth{s}}

     \newcommand{\C}{\Bbb{C}}

     \newcommand{\Q}{\Bbb{Q}}
     
     \newcommand{\Z}{\Bbb{Z}}

    \newcommand{\qed}{{\hfill$\Box$}}

    \newcommand{\ul}[1]{\underline{#1}}
    
    \newcommand{\ti}[1]{\tilde{#1}}

    \newcommand{\me}{^{-1}}
    \newcommand{\mal}{^{\times}}
    
    \newcommand{\df}{\stackrel{\mathrm{def}}{=}}
    \newcommand{\mr}{\mathrm}
    \newcommand{\clo}{^{\mr{c}}}
    \newcommand{\ilim}[1]{\raisebox{-1mm}{$\lim\atop{\leftarrow\atop\raisebox{0.5mm}{{\tiny $#1$}}}$}}

    \newcommand{\zp}{{\Bbb{Z}_p}}
    \newcommand{\zl}{{\Bbb{Z}_l}}
    \newcommand{\qp}{{\Bbb{Q}_p}}

    \newcommand{\lto}{\longrightarrow}
    \newcommand{\da}{\downarrow}

    \newcommand{\pht}{\phantom}

    \def\daz#1{#1\da\pht{#1}}
    
    \newcommand{\ga}{\gamma}
    \newcommand{\Ga}{\Gamma}
    \newcommand{\si}{\sigma}
    \newcommand{\Si}{\Sigma}
    \newcommand{\de}{\delta}
    \newcommand{\De}{\Delta}
    \newcommand{\ze}{\zeta}
    \newcommand{\la}{\lambda}
    \newcommand{\La}{\Lambda}
    \newcommand{\be}{\beta}

    \newcommand{\al}{\alpha}
    \newcommand{\ve}{\varepsilon}

    \newcommand{\noi}{\par\noindent}
    \newcommand{\sn}{\par\smallskip\noindent}
    \newcommand{\mn}{\par\medskip\noindent}
    \newcommand{\bn}{\par\bigskip\noindent}
    \newcommand{\bbn}{\par\bigskip\bigskip\noindent}
    
    \newcommand{\Section}[2]{\bbn {\large #1\,. \ {\sc #2}}
                             \nopagebreak
                             \nz}

    \newcommand{\nf}[2]{\\[1.5ex]
                        \bmp{1cm}
                         (#1)
                        \emp % \hsp{30}
                        \bmp{13.5cm}
                         \bct
                          $#2$
                         \ect
                        \emp\\[1.5ex]
         }

    \newcommand{\nz}{\\[1ex]}

    \newcommand{\hsp}[1]{\hspace*{#1mm}}

    \newcommand{\mmargin}{
     \textheight 230truemm
     \textwidth 155truemm
     \topmargin -10truemm
     \oddsidemargin 5truemm
     \evensidemargin 5truemm
     }

    \newcommand{\bmp}{\begin{minipage}}
    \newcommand{\emp}{\end{minipage}}
    \newcommand{\btb}{\begin{tabular}}
    \newcommand{\etb}{\end{tabular}}
    \newcommand{\barr}{\begin{array}}
    \newcommand{\earr}{\end{array}}
    \newcommand{\bit}{\begin{itemize}}
    \newcommand{\eit}{\end{itemize}}
    \newcommand{\ben}{\begin{enumerate}}
    \newcommand{\een}{\end{enumerate}}
    \newcommand{\bct}{\begin{center}}
    \newcommand{\ect}{\end{center}}
    \newcommand{\bfr}{\begin{flushright}}
    \newcommand{\efr}{\end{flushright}}
    \newcommand{\bea}{\begin{eqnarray*}}
    \newcommand{\eea}{\end{eqnarray*}}
    \newcommand{\bqo}{\begin{quote}}
    \newcommand{\eqo}{\end{quote}}
    \newcommand{\bdc}{\begin{description}}
    \newcommand{\edc}{\end{description}}
    \newcommand{\bdia}{\begin{CD}}
    \newcommand{\edia}{\end{CD}}

    \definecolor{light}{gray}{.3}

    \newcommand{\tr}{\mathrm{tr\,}}

    \newcommand{\ind}{\mathrm{ind\,}}
    \newcommand{\infl}{\mathrm{infl\,}}

    \newcommand{\res}{\mathrm{res\,}}
    
    \newcommand{\Det}{\mathrm{Det\,}}
    \newcommand{\sr}[2]{{\,\stackrel{#1}{#2}\,}}
    \newcommand{\fra}[2]{{\,\frac{#1}{#2}\,}}

    \newcommand{\theorem}{\sn
           \bdc
           \item[{\sc Theorem.}] \em }
    \newcommand{\Theorem}[1]{\sn
           \bdc
           \item[{\sc Theorem {#1}.}]  \em }

    \newcommand{\Stop}{\edc \sn\rm} % GROSSES S !!

    \newcommand{\Lemma}[1]{\sn
           \bdc
           \item[{\sc Lemma {#1}.}] \em }
    
    \newcommand{\Proposition}[1]{\sn
           \bdc
           \item[{\sc Proposition {#1}.}] \em }

    \newcommand{\proof}{{\sc Proof.} \ }
    \newcommand{\remark}{\sn{\sc Remark.} \ }

\mmargin

 \def\zpp{{\Z_p}}
 \def\zpgs{{\zpp[[G_S]]}}
 \def\zphs{{\zpp[[H_S]]}}
 \def\zphsp{{\zpp[[H_S^+]]}}
 \def\qq{{\mathcal{Q}}}
 \def\ab{{\mr{ab}}}
 \def\ver{{\mr{ver}}}
 \def\tr{{\mr{tr}}}
 
 \def\lw{{\La_\wedge}}

\def\zefe{{\ti\ze_F(1-k,\ve)}}
\def\zefeg{{\ti\ze_F(1-k,\ve_g)}}
\def\cn{{{\cal{N}}}}
\def\dex{{\de^{(x)}}}

\def\muu{{m(U)}}
\def\pmu{{p^\muu}}

\def\cnf{{\cn_F}}
\def\cnk{{\cn_K}}
\def\sif{{\Si_F}}
\def\mus{{\mu_\Si}}
\def\gfs{{G(F_S/F)}}
\def\osf{{|\Si_F|}}
\def\osfk{{|\Si_F|k}}
\def\verfl{{\ver_F^L}}
\def\verkf{{\ver_K^F}}

\def\ee{{\mathcal{E}}}
 \def\mm{{\mathcal{M}}}
 \def\gann{{\Ga_{00}}}
 
 \def\el{{\ve_L}}
 \def\stab{{\mr{St}}}
 \def\zkp{{\Z_{(p)}}}
\def\pu{{p^{m_F(U)}}}
\def\gf{{G(F_S/F)}}

\begin{document}

\title{Congruences between abelian pseudomeasures\,, II}
 \author{Jürgen Ritter \ $\cdot$ \ Alfred Weiss \
 \thanks{We acknowledge financial support provided by DFG and NSERC.}
 }
\date{\pht\today}

 \maketitle
 \pagestyle{myheadings}
\markright{Congruences between abelian pseudomeasures}
% \bct\bmp{14.5cm}{\small We extend the main result of [Math. Res. Lett. 15
% (2008), 715-725] to Galois extensions $L/K$ of totally real number
% fields of arbitrary prime power degree, thereby providing a basic
% tool for our proof of the `main conjecture' of equivariant Iwasawa
% theory [RW10]. \mn 2000 Mathematics Subject Classification: 11R23,
% 11R42}\emp\ect

 \bn In this paper we generalize the main result of [RW7]. As
 there, $K$ is a totally real number field (finite over $\Q$), $p$
 a fixed prime number, and $S$ a fixed finite set of non-archimedean
 primes of $K$ containing all primes above $p$. Let $K_S$ denote
 the maximal abelian extension of $K$ which is unramified (at all
 non-archimedean primes) outside $S$ and set $G_S=G(K_S/K)$.
 Serre's pseudomeasure $\la_K=\la_{K,S}$ has the property that
 $(1-g)\la_K$ is in the completed group ring $\zpgs$ for all $g\in
 G_S$ [Se].
 \mn Let $L$ be a totally real Galois extension of $K$ of
 $p$-power degree with group $\Si=G(L/K)$. We require that
 the finite set $S$ contains all the primes of $K$ which ramify in
 $L$. Letting $F$ run through the intermediate fields of $L/K$, we
 denote by
 \ben\item[] $F_S$ the maximal abelian extension of $F$
 unramified outside the primes of $F$ above $S$\,,
 \item[] $F_S^+$ its maximal totally real subfield (hence
 $F\subset F_S^+$)\,,
 \item[] $\la_{F,S}$ Serre's pseudomeasure with respect to $F$ and
 $S$\,,
 \item[] $H_S=G(L_S/L)\,,\,H_S^+=G(L_S^+/L)$\,.\een
 Note that $G(L_S/F)^\ab=G(F_S/F)$ and, for $F\subseteq F'$, that $G(L_S/F')$ is an open subgroup of
 $G(L_S/F).$ This yields the transfer map $G(F_S/F)\to G(F_S'/F')$\,, and in particular,
 $$\ver_K^F:G_S\to G(F_S/F)\,,\,
 \ver_F^L:G(F_S/F)\to H_S\,,\,\ver_K^L=\ver_F^L\circ\ver_K^F:G_S\to H_S\,.$$
 In order to state the new main result we still need to recall the
 definition of the Möbius function $\mu=\mu_\Si$ of the poset of
 subgroups of the finite group $\Si$. It is defined by
 $$\mu(1)=1\,,\,\mu(\Si')=-\sum_{1\le\Si''<\Si'}\mu(\Si'') \
 \mr{for} \ 1\neq\Si'\le \Si\,.$$
 For $K\subseteq F\subseteq L$ and $g\in G_S$ set
 $\ti\la_F=2^{-[F:\Q]}\la_{F,S}$ and
 $\ti\la_{g_F}=(1-g_F)\ti\la_F$\,, where $g_F=\ver_K^F g$.
 Moreover, denote the Galois group of the cyclotomic
 $\zp$-extension $F_\infty/F$ by $\Ga_F$.
 \theorem There exist $g\in G_S$ so that $g_F$ has image $\neq1$
 in $\Ga_F$ for $K\subseteq F\subseteq L$ and the image of
 $$\sum_{K\subseteq F\subseteq
 L}\mus(G(L/F))\ver_F^L(\ti\la_{g_F})\,,$$
 under $\zphs\to\zphsp$\,, is in the trace ideal $\tr_\Si(\zphsp$
 of $\zphsp^\Si$\,.\Stop
 The theorem and its proof generalize [RW7], where $\Si$ had order
 $p$, and therefore are an application of the methods of [DR] in
 the language of [Se]. The precise formulation,
 however, is dictated by the needs of equivariant Iwasawa theory
 [RW2,p.564].
 \mn More precisely, there one
 has a profinite extension $\mathbf{K}/K$ of totally real fields with
 Galois group $G$, with $K/\Q$ finite and $\mathbf{K}$ containing the cyclotomic
 $\zpp$-extension $K_\infty$ of $K$, such that\linebreak
 $[\mathbf{K}:K_\infty]<\infty$\,.
 In this situation, the so-called {\em `main conjecture' of
 equivariant Iwasawa theory}~\footnote{see also e.g.~[CFKSV; FK; Ha1,2; K; Ka; RW3,7,8]} concerns the {\em Iwasawa
 $L$-function}
 $L_{\mathbf{K}/K,S\cup S_\infty}$\,, which is built from
 the $S\cup S_\infty$-truncated \footnote{$S_\infty$ is the set of archimedean primes of $k$} $p$-adic Artin $L$-functions
 $L_{p,K,S\cup S_\infty}(s,\chi)$ attached to the $\Q_p\clo$-irreducible
 characters $\chi$ of $G$ with open kernel (see loc.cit.~-- but
 with \,$\mathbf{K}\,,\,K\,,\,p$ \,denoted \,$K\,,\,k\,,\,l$\,, respectively).
 Assuming, for odd $p$, Iwasawa's conjecture that the $\mu$-invariant for $\mathbf{K}/K$ vanishes,
 the `main conjecture', up to its uniqueness assertion, is equivalent to
 $$L_{\mathbf{K}/K,S\cup S_\infty}\quad\mr{is\ in}\quad\Det(K_1(\lw G))\,,$$
 where $\lw G$ is the
 completion of the localization of $\zpp[[G]]$ obtained by
 inverting all elements which are regular modulo $p$
 (see [RW3, Theorem A]
 \footnote{for $p=2$, compare [RoW]}).

 \sn Now assuming that $G$ is a
 $p$-elementary group, i.e., a direct product of a finite cyclic group of order prime to $p$ and a pro-$p$ group, and picking
 an abelian normal open
 subgroup $A$ of $G$ with factor
 group $\Si=G/A$ a $p$-group,
 the {\em Möbius-Wall congruence}
 \nf{MW}{\sum_{A\le U\le
 G}\mu_\Si(U/A)\ver_U^A(\res_G^Uy)\equiv0\mod\tr_\Si(\lw A)}
 %is valid
 holds for all units $y\in\lw G$.  (MW)
 is in [RW10, \S2].
 \sn If the `main conjecture' were true, so  Det$(y)=L_{\mathbf{K}/K,S\cup S_\infty}$
 with a $y\in K_1(\lw G)$, then on plugging $y$ into (MW) we would obtain $$\sum_{A\le U\le
 G}\mu_\Si(U/A)\ver_{U^\ab}^A(\la_{U^\ab})\equiv0\mod\tr_\Si(\lw
 A)\,.$$ And this {\em abelian pseudomeasure congruence} turns
 out to be a consequence of the Theorem stated above. This is addressed in
 detail in [RW10], where the `main conjecture' for arbitrary extensions $\textbf{K}/K$, up to its uniqueness assertion, is verified.

 \mn The organization of the paper parallels [RW7] to the extent
 that, with two exceptions \footnote{Propositions 5 and 9}, the
 numbering and, mutatis mutandis, statements of the lemmas and propositions
 persist, with proofs omitted if they are essentially already in
 [RW7]. The new ingredient is the identification of the congruence
 $(4\star)$ of Proposition 4 as the difference of constant terms of
 $q$-expansions at two cusps of a Hilbert modular form of
 Eisenstein type. This modular form is exhibited in \S2 and then
 studied via the $q$-expansion principle of Deligne and Ribet; the
 hypothesis of Proposition 4 is deduced, in Proposition 9, from a
 property of Möbius coefficients in [HIÖ]. The proof, in \S3, of
 the main result is then a computation of constant terms of
 $q$-expansions at the ``special cusps'' of Proposition 5.
 \sn Our special cusps are a simple device to avoid comparing
 constant coefficients of $q$-expansions of $F$ and $F_{|k}U_\be$
 at {\em arbitrary} cusps in Lemma 6. Having overlooked the need for
 this comparison in [RW7] implies that we now have its Theorem only for special $g=g_K\in G_S$.
 However, the present theorem is better suited for equivariant Iwasawa theory;
  see the
 Remark in \S3.
 \Section{1}{A sufficient condition for a pseudomeasure congruence}
 For a coset $x$ of an open subgroup $U$ of $G(F_S/F)$ set
 $\dex(g)=1$ or 0 according as $g\in x$ or not. Then, for even integers
 $k\ge1$, define \ $\ti\ze_F(1-k,\dex)=2^{-[F:\Q]}\ze_{F,S}(1-k,\dex)\in\Q$ \ to
 be $2^{-[F:\Q]}$ times the value at $1-k$ of the partial $\ze$-function for the set
 of integral ideals $\fa$ of $F$ prime to $S$ with Artin symbol
 $(\fa,F_S/F)$ in $x$. Note that the
 definition of $\ti\ze_F(1-k,\de^{(x)})$ extends linearly to locally constant functions $\ve$
 on $G(F_S/F)$ with values in a $\Q$-vector space and gives values
 $\ti\ze_F(1-k,\ve)$ in that vector space, as usual.
 \mn Let $\cn=\cn_{F,p}:G(F_S/F)\to\zp\mal$ be that continuous character whose
 value on $(\fa,F_S/F)$ for an integral ideal $\fa$ of $F$ prime to $S$ is its
 absolute norm $\cnf\fa$.
 For $g\in G(F_S/F)\,,\,k\ge1$ and $\ve$ a locally constant
 $\qp$-valued function on $G(F_S/F)$ we define, following [DR],
 $$\ti\De_g(1-k,\ve)=\zefe-\cn(g)^k\zefeg\in\qp\,,$$
 where $\ve_g(g')=\ve(gg')$ for $g'\in G(F_S/F)$.
 \Theorem{{\rm[(0.4) of [DR]]}} Let $\ve_1,\ve_2,\ldots$ be a
 finite sequence of locally constant functions \linebreak $G(F_S/F)\to\qp$ so that
 \ $\sum_{k\ge1}\ve_k(g')(\cn g')^{k-1}\in\zp$ for all $g'\in G(F_S/F)$\,. Then
  $$\sum_{k\ge1}\ti\De_g(1-k,\ve_k)\in\zp\quad\mr{for\ all}\ g\in G(F_S/F)\,.\,\footnote{For $p=2$ compare [RW7,\S5]}$$\Stop
 Call an open subgroup $U$ of $G(F_S/F)$ {\em admissible}, if $\cn(U)\subset
 1+p\zp$, and define $m_F(U)\ge1$ by
 $\cn(U)=1+p^{m_F(U)}\zp$.
 \Lemma{1} If $U$ runs through the cofinal system of admissible open subgroups
 of $G(F_S/F)$, then $\zp[[G(F_S/F)]]=\ilim{U}\zp[G(F_S/F)/U]/p^{m_F(U)}\,\zp[G(F_S/F)/U]$\,.\Stop
% See [RW7,Lemma 1].
 \Proposition{2} For $h\in G(F_S/F)$ there is a unique element
 $\ti\la_{h}\in\zp[[G(F_S/F)]]$, independent of $k$, whose image in
 $\zp[G(F_S/F)/U]/p^{m_F(U)}$ is $$\sum_{x\in G(F_S/F)/U}\ti\De_{h}(1-k,\dex)\cn
 (x)^{-k} x\mod p^{m_F(U)}\zp[G(F_S/F)/U]$$ for all admissible $U$, where $\mathcal{N}$ here also denotes
 the homomorphism $G(F_S/F)/U\to(\zp/\pmu)\mal$ induced by our previous $\mathcal{N}$. Moreover, if $\ti\la_F$ is $2^{-[F:\Q]}$ times the pseudomeasure of
 [Se], then $(1-h)\ti\la_F=\ti\la_{h}$\,. \Stop
% The proposition follows from the theorem above, see [RW7, Proposition 2].
 \Lemma{3} \ben\item[(1)]Let $V$ be an admissible open subgroup of $H_S$. If $U\le(\ver_F^L)\me(V)$\,, then  $m_F(U)\ge m_L(V)-e_F$
 where $[L:F]=p^{e_F}$\,.
 \item[(2)]  Let $s\in G(L_S/K)$ be an extension of $\si\in \Si$. Then $(F^\si)_S=(F_S)^{s}$ and, setting $g_F=\verkf(g)$ for $g\in G_S$\,,
 \ $g_F^s=g_{F^\si}$\,. Moreover,
 $\ver_F^L(\ti\la_{g_F})^{s}=\ver_{F^\si}^L(\ti\la_{g_{F^\si}})$\,.
 \een\Stop
 Statement (1) is due to $\cn_L(\ver_F^L
 g)=\cn_F(g)^{[L:F]}$ for $g\in G(F_S/F)$ \,.
 For (2), the first claim follows from $(F_S)^{s}\supseteq
 F^{s}=F^\si$ and $G((F_S)^{s}/F^{s})=
 G(F_S/F)^{s}$\,, the latter implying
 $(F_S)^{s}\subseteq(F^\si)_S$ and then
 $F_S\subseteq(F^\si)_S^{s\me}\subseteq F_S$\,, hence equality
 everywhere.
 \sn The second claim is a direct consequence of the definition of
 group transfer `ver'. Namely $\verkf(g)$ is a certain product built with respect
 to coset representatives of $G(L_S/F)$ in $G(L_S/K)$\,, which
 $s$ takes to coset representatives of $G(L_S/F^\si)$ in
 $G(L_S/K)$\,, whence the multiplicativity of $s$ yields
 $\verkf(g)^{s}=\ver_K^{F^\si}(g^{s})$\,. But $g\in
 G_S=G(L_S/K)^\ab$ implies $g^{s}=g$.
 \sn Finally, concerning the last claimed equality, by Proposition 2 and $\cn_{F^\si}(x^s)=\cnf(x)$ it suffices to show
 $$\ti\De_{g_{F^\si}}(1-k,\de_{F^\si}^{(x^s)})=\ti\De_{g_F}(1-k,\de_F^{(x)})$$
 for $x$ a coset of any admissible open $G(L/F)$-stable subgroup $U$ of $G(F_S/F)$\,.
 Thus it suffices to show $$\ti\ze_{F^\si}(1-k,\de_{F^\si}^{(x^s)})=\ti\ze_F(1-k,\de_F^{(x)})$$
 for all such $x$, since\,\footnote{recall our convention $\ve_g(g')=\ve(gg')$}
 \ $(\de_{F^\si}^{(x^s)})_{g_{F^\si}}=\de_{F^\si}^{(g_{F^\si}\me x^s)}=\de_{F^\si}^{((g_{F}\me x)^s)}$
 and $(\de_F^{(x)})_{g_F}=\de_F^{(g_F\me x)}$\,. Viewing $\de_F^{(x)}$ as a complex valued
 function on $\gfs/U$ and writing it as a $\C$-linear combination of
 the (abelian) characters $\chi$ of $\gfs/U$, it now suffices
 to check that $\ti\ze_F(1-k,\chi)=\ti\ze_{F^\si}(1-k,\chi^s)$\,,
 and this follows from the
 compatibility of the Artin $L$-functions with inflation and induction. Indeed,
 inflating $\chi$ from $G(F_S/F)/U$ to $G(F_S/F)$ and further to $G(L_S/F)$ and then inducing up to
 $G(L_S/K)$, and analogously with $\chi,F,U$ replaced by
 $\chi^s,F^\si,U^s$ (note that $U^s$ is well-defined), we have
 $$\ind_{G(L_S/F)}^{G(L_S/K)}\infl_{G(F_S/F)/U}^{G(L_S/F)}(\chi)=\ind_{G(L_S/F^\si)}^{G(L_S/K)}\infl_{G(F_S^s/F^\si)/U}^{G(L_S/F^\si)}(\chi^s)\ .$$
 \sn Lemma 3 is established.\qed
 \bn Set $\sif=G(L/F)$\,, and,
 for any set $X$ carrying a natural $\Si$-action, denote
 the stabilizer subgroup
 of $x\in X$ in $\Si$ by $\stab_\Si(x)=\{\si\in\Si:\si(x)=x\}$\,.
 Also, $\zkp\subset\Q$ is the localization of $\Z$ at its
 prime ideal $p\Z$.
 \Proposition{4} For $g\in G_S$ set $g_F=\verkf g$ for all
 $K\subseteq F\subseteq L$\,. Then
 \mn
 $(4\star)\hsp{15}\sum\limits_{K\subseteq F\subseteq L}\mus(\Si_F)\ti\De_{g_F}(1-|\Si_F|k,\ve_L\ver_F^L)\equiv0\mod|\stab_\Si(\ve_L)|\zkp\,,$
 \mn  for all even locally constant $\Z_{(p)}$-valued functions $\ve_L$ on $H_S$\,, implies that
 $$\fs_g\df\sum_F\mus(\Si_F)\ver_F^L(\ti\la_{g_F})$$
 has image, under the map $\zp[[H_S]]\to\zp[[H_S^+]]$, in $\tr_\Si(\zp[[H_S^+]]$\,.
 \Stop
 For the proof of the proposition we first recall that a locally
 constant function $\el$ on $H_S$ is {\em even} if
 $\el(c_wh)=\el(h)$ for all $h\in H_S$ and all `Frobenius
 elements' $c_w$ at the archimedean primes $w$ of $L$, i.e., at the
 restrictions $c_w\in H_S$ of complex conjugation with respect to
 the embeddings $L_S\hookrightarrow\C$ inducing $w$ on $L$. We
 denote by $C$ the group generated by the $c_w$'s, so $H_S^+=H_S/C$\,.
 \sn We next observe that $\fs_g\in\zphs^\Si$\,.
 Namely,
 $\ver_F^L(\ti\la_{g_F})^s=\ver_{F^\si}^L(\ti\la_{g_F}^s)=\ver_{F^\si}^L(\ti\la_{g_{F^\si}})$\,,
 by (2) of Lemma 3. Moreover, $\mus(\sif)=\mus(\Si_{F^\si})$\,.
 \sn Turning finally to the image of $\fs_g$ under the map
 $\zp[[H_S]]\to\zp[[H_S^+]]$\,, we first replace the diagram in
 [RW7,p.718] by the diagram below, in which $N=\ker\ver_F^L$\ :
 $$\barr{cccc}\zp[[G(F_S/F]]&\to&\ilim{U\ge N}\Z_p[G(F_S/F)/U]/p^{m_F(U)}&\\
 \daz{\ver_F^L}&&\da&\\
 \zp[[H_S]]&\sr{\simeq}{\to}&\ilim{\Si-\mr{stable}\,V}\Z_p[H_S/V]/p^{m_L(V)-e_F}&.\earr$$
 Recall here that the right vertical map takes $(x_U)_U$ to
 $(y_V)_V$ by means of
 $$\Xi\,:\quad\zp[G(F_S/F)/U]/p^{m_F(U)}\sr{\ver}{\lto}\zp[H_S/V]/\pu\to\zp[H_S/V]/p^{m_L(V)-e_F}\,,$$
 whenever $U\le(\verfl)\me(V)$. \footnote{Note that 'ver' is the
 $\zp$-linear map induced by the group homomorphism obtained by
 factoring $\gf\sr{\verfl}{\to}H_S\to H_S/V$ through
 $\gf\to\gf/U$\,.}
 \mn Since the $m_L(V)$'s are unbounded, there are admissible open
 $\Si$-stable $V\le H_S$ with $m_L(V)-e_F\ge e_K\ (\,\forall\,
 F\,)$. For any such $V$ then $\Z_p[H_S/V]/p^{m_L(V)-e_F}$ maps
 onto $\Z_p[H_S/V]/|\Si|$ and we write the image of $\fs_g$ in
 here as \ $\sum_{y\in H_S/V}c_yy$\,. Because $\Si$ fixes $\fs_g$\,,
 $c_{y^\si}=c_y$ for all $\si$. Since
 $\sum_{\si\in\Si\,\mr{mod}\,\stab_\Si(y)}c_{y^\si}y^\si=c_y\sum
 y^\si$\,, it follows that $\fs_g$
 will be in $\tr_\Si(\Z_p[H_S/V])+|\Si|\Z_p[H_S/V]$ provided that
 $$c_y\equiv0\mod|\stab_\Si(y)|\ .$$
 We compute the coefficient $c_y$ of the image of $\fs_g$. By
 Proposition 2 with $U=(\ver_F^L)\me(V)$ (compare [RW7, (ii) on p.719] noting that any $k$ is allowed)
 \ $\ver_F^L(\ti\la_{g_F})$ has image
 $$\sum_{x\in G(F_S/F)/U}\ti\De_{g_F}(1-|\Si_F|k,\de_F^{(x)})\cn_F(x)^{-|\Si_F|k}\ver_F^L(x)\,,$$
 where $\de_F^{(x)}$ is the characteristic function of the coset
 $x\subseteq G(F_S/F)$\,.
 \sn Since $G(F_S/F)\sr{\ver_F^L}{\to}H_S\to H_S/V$ has kernel
 $U=(\ver_F^L)\me(V)$, either $y$ is not in the image of $\ver_F^L$
 or $y=\ver_F^L(x^{(F)})$ for a unique $x^{(F)}\in G(F_S/F)/U$\,.
 Note that $\de_L^{(y)}\ver_F^L$ is $=0$ in the first case and
 $=\de_F^{(x^{(F)})}$ in the second, when also
 $\cn_F(x^{(F)})^{-|\Si_F|k}=\cn_L(\ver_F^L (x^{(F)}))^{-k}$\,. Thus, Möbius-summing over $F$ we obtain
 $$c_y=(\,\sum_F\mus(\sif)\ti\De_{g_F}(1-|\Si_F|k,\de_L^{(y)}\ver_F^L)\,)\cn_L(y)^{-k}\,.$$
 Our hypothesis $(4\star)$ now implies that $\fs_g$ is in
 $\tr_\Si(\Z_p[H_S/V])+|\Si|\Z_p[H_S/V]$ for all $V\ge C$ (recall that
 $\de_L^{(y)}$ is even when $V\ge C$ and $y\in H_S/V$). Since $\fs_g$ is fixed by $\Si$
 and $|\Si|\zp[H_S/V]^\Si\subseteq\tr_\Si(\zp[H_S/V])$\,, it follows that $\fs_g\in\tr_\Si\Z_p[H_S/V]$\,.
 \sn This finishes the proof of Proposition 4.\qed
  \Section{2}{Applying the $q$-expansion principle of [DR]}
 Given an even integer $k$ and an even locally constant
 $\zkp$-valued function $\el$ on $H_S$, choose an open subgroup $V$ of $H_S$ so that $\el$ is
 constant on each coset $H_S/V$ and let
 $\ff\subset|\Si|\fo_K$ be an integral ideal, with all its
 prime factors contained in $S$, so that, for all
 $K\subseteq F\subseteq L$, $\ff\fo_F$ is a multiple of the conductor of the
 field fixed by $(\verfl)\me(V)$ acting on $F_S$.
 \mn As in [DR,p.229] we write $\hat K$ for the ring of `finite' ad\`eles
 of $K$ and let $j:\hat K\mal\to G_S\,,\psi:\hat K\mal\to\ilim{\ff'}G_{\ff'}$\,, with
 $\ff'$ running over the multiples of $\ff$ that have all their prime divisors in $S$, be the maps defined in [DR,p.243].
 \Proposition{5} There exist $\ga\in\hat K\mal$ so that
 \ben\item[a)] $\ga$ and $\ga\me$ are in $1+\hat\ff$\,, and
 \item[b)] the image $g\in G_S$ of $\ga$ under $j:\hat K\mal\to
 G_S$ has $\verkf(g)\in G(F_S/F)$ \ul{not} in the kernel of
 $G(F_S/F)\to \Ga_F$\,, for all $K\subseteq F\subseteq L$.\een\Stop
 \proof Choose $f>0$ in $\Z\cap\ff$, and let $\al\in\hat\Q\mal$ be
 the `finite' id\`ele with components $1+f$ (respectively 1) at
 primes $q|f$ (respectively $q\nmid f$). Then $\al$ and $\al\me$
 belong to $1+\widehat{\Z f}$. For every extension
 $\Q\hookrightarrow F$ let $\al_F$ be the image of $\al$ under the diagonal
 inclusion $\hat\Q\mal\hookrightarrow\hat F\mal$.  We use $\ilim{\ff'}G_{\ff'}=G_S$, on identifying the
 inverse limit of ray class groups with the one of the corresponding Galois
 groups, as is conventional in [DR,p.240]; so $j=\psi=(\psi_{\ff'})_{\ff'}$ by [DR,2.33].
 \mn Now $\psi_{\ff'}(\al_K)=\psi_{\ff'}((1+f)\me\al_K)$\,, since
 $1+f\in K^{\gg0}$\,, and $(1+f)\me\al_K$ has $\fq$-component 1
 (respectively $(1+f)\me$) when $\fq|\ff$ (respectively
 $\fq\nmid\ff$), so that the (fractional) ideal of $\fo=\fo_K$
 `generated' by $(1+f)\me\al_K$ is $\fa=(1+f)\me\fo$\,, which is coprime to
 $\ff$.
 \sn Recall that $K_S$ contains the $\zp$-cyclotomic extension
 $K_\infty$
 of $K$. By [Se, 2.2], $g=\psi(\al_K)\in G_S$ acts on $p$-power
 roots of unity as $g=(\fa,K_S/K)$, i.e., by raising them to the
 power $\cn(\fa)=(1+f)^{-[K:\Q]}$\,. It follows that the image of
 $g$ under $G_S\to\Ga_K$ acts non-trivially, so {\em a)} holds for $\al_K$ and also {\em b)} when $K=F$.
 \sn The argument for $F$ is the same with $K,\ff,S$
 replaced by $F,\ff\fo_F,S_F$ and shows that $\psi_F(\al_F)$ acts on
 $p$-power roots of unity by $(1+f)^{-[F:\Q]}$. Note also that
 $\psi_F(\al_F)=\verkf(g)$ by the usual relation between inclusion and transfer. Setting
 $\ga=\al_K$ we then get {\em a)} and {\em b)}.
 \sn This completes the proof of Proposition 5.\qed
 \mn The next three results of [RW7] concern Hilbert modular forms
 with emphasis on their $q$-expansions. Thus Lemma 6 constructs a
 Hecke operator $U_\be$ on $\mm_k(\gann(\ff),\C)$\,, Lemma 7
 discusses restriction
 $\res_F^K:\mm_k(\gann(\ff\fo_F),\C)\to\mm_{[F:K]k}(\gann(\ff),\C)$
 for field extensions $F/K$, and Proposition 8 is our bridge to
 [DR].
 \bn With $k,\ve_L$ and $\ff$ as at the beginning of the section and any $g\in G_S$, we next exhibit a Hilbert modular
 form $\ee$ in $\mm_{|\Si|k}(\gann(\ff),\C)$ with the constant term of its standard
 $q$-expansion equal to \
 $\sum_F\mus(\Si_F)\ti\De_{g_F}(1-|\Si_F|k,\ve_L\ver_F^L)$ (compare $(4\star)$).
 \mn First, by [DR, (6.1)], in the form of
 Proposition 8, there are modular forms
 $$E_F\df G_{\osfk,\el\verfl}\in\mm_\osfk(\gann(\ff\fo_F),\C)$$
 of weight $\osfk$ with standard $q$-expansion
 $$\ti\ze_F(1-\osfk,\el\verfl)+\sum_{\substack{\mu\gg0\\
 \mu\in\fo_F}}\Big(\sum_{\substack{\mu\in\fa\subseteq\fo_F\\
 \fa\,\mr{prime\,to}\,S}}\el\verfl(\fa)\cnf(\fa)^{\osfk-1}\Big)\,q_F^\mu\ .$$
 Appealing to Lemma 7 and Lemma 6, we
 apply $\res_F^K$ and the Hecke operator $U_{[F:K]}$ to the
 modular form $E_F$ displayed above, and obtain, for each $F$, the new
 modular form $$\ee_F=(\res_F^K E_F)_{\mid_{[F:K]|\osfk}}U_{[F:K]}$$ of weight $|\Si|k$ in
 $\mm_{|\Si|k}(\gann(\ff),\C)$ and with standard $q$-expansion
 $$\ti\ze_F(1-\osfk,\el\verfl)+\sum_{\substack{\al\gg0\\\al\in\fo_K}}\Big(\sum_{[\al]_F}\el(\fa_F\fo_L)\cnf(\fa_F)^{[L:F]k-1}\Big)q_K^\al\,,$$
 where, for $K\subseteq F\subseteq L$, $[\al]_F$ denotes the set of all pairs $(\al_F,\fa_F)$
 satisfying $$0\ll\al_F\in\fa_F\subseteq\fo_F\,,\,\fa_F\ \mr{prime\
 to}\ S\,,\,\tr_{F/K}(\al_F)=[F:K]\al\ .$$
 Here we have used $(\el\circ\verfl)(\fa_F)=(\el\verfl)(\fa_F,F_S/F)=\el(\fa_F\fo_L)$.
 \mn We Möbius-sum all these $\ee_F$ and arrive at the modular form
 $$\ee\df\sum_{K\subseteq F\subseteq L}\mus(\sif)\ee_F$$
 in $\mm_{|\Si|k}(\gann(\ff),\C)$ whose standard $q$-expansion has constant coefficient
 \nf{0}{\sum_F\mus(\sif)\ti\ze_F(1-\osfk,\el\verfl)}
 and higher coefficients
 \nf{$\al$}{\sum_F\mus(\sif)\sum_{[\al]_F}\el(\fa_F\fo_L)\cnf(\fa_F)^{[L:F]k-1}\,,}
 at $0\ll\al\in\fo_K$\,.
 \Proposition{9} Assume that $k$ is an even positive integer and
 $\el$ an even locally
 constant $\zpp$-valued function on $H_S$. Then $$\sum_{K\subseteq F\subseteq L}\mus(\sif)\sum_{[\al]_F}\el(\fa_F\fo_L)\cnf(\fa_F)^{[L:F]k-1}
 \equiv0\mod|\stab_\Si(\el)|\Z_{(p)}\ .$$\Stop
 \proof  Utilizing the natural action of $\Si$ on the set
 $$[\al]_L=\{(\al_L,\fa_L):0\ll\al_L\in\fa_L\subseteq\fo_L\,,\,\fa_L\
 \mr{prime\ to}\ S\,,\,\tr_{L/K}(\al_L)=[L:K]\al\}\,,$$
 given by \ $(\al_L,\fa_L)^\si=(\al_L^\si,\fa_L^\si)$\,, we
 identify the set
  $[\al]_F$ with the subset $[\al]_L^\sif$ of $\sif$-fixed
 points in $[\al]_L$ by means of the map $$\iota_F:[\al]_F\to[\al]_L\,,\quad(\al_F,\fa_F)\mapsto(\al_F,\fa_F\fo_L)\ .$$
 Indeed, $\iota_F$ is obviously injective
 and has image $[\al]_L^\sif$ because
 \ben\item[] $\al_L\in F$ if $\al_L^\si=\al_L$ for all
 $\si\in\sif$\,,
 \item[] $\fa_L=\fa_L^\si$ for all $\si\in\sif$ implies
 $\fa_L=\fa_F\fo_L$ for some integral ideal $\fa_F$ of $F$, since
 $\fa_L$ is prime to $S$ and whence every prime divisor of $\fa_L$
 is unramified in $L/K$\,.
 \een As a first consequence, formula $(\al)$ can be rewritten as
 \nf{$\al'$}{\sum\limits_{K\subseteq F\subseteq
 L}\mus(\sif)\sum\limits_{(\be,\fb)\in[\al]_L^\sif}\el(\fb)\cn_L(\fb)^{k-\fra1{\osf}}\,,}
 as $\cn_L(\fa_F\fo_L)=\cnf(\fa_F)^{[L:F]}$\,.
 \mn We isolate the part of the sum $(\al')$ that belongs to a
 fixed $(\be,\fb)\in[\al]_L$\,. It is
 $$\sum_{F\,\mr{so}\,\sif\le\stab_\Si(\be,\fb)}\mus(\sif)\el(\fb)\cn_L(\fb)^{k-\fra1{\osf}}\,,$$
 and correspondingly, with $(\be,\fb)$ replaced by
 $(\be,\fb)^\si$\,,
 $$\sum_{F\,\mr{so}\,\sif\le\stab_\Si(\be,\fb)^\si}\mus(\sif)\el(\fb^\si)\cn_L(\fb^\si)^{k-\fra1{\osf}}\,.$$
 The group $\Si$ acts on $H_S$ and so on $\el$ by
 $\el^\si(h)=\el(h^{\si\me})$\,. We now consider the part of the
 sum $(\al')$ that belongs to the
 $\stab_\Si(\el)$-orbit of $(\be,\fb)$\,:
 \nf{i}{\sum\limits_{\si\in[\stab_\Si(\el)\cap\stab_\Si(\be,\fb)\setminus\stab_\Si(\el)]}
 \el(\fb)\cn_L(\fb)^k\sum\limits_{F\,\mr{so}\,\sif\le\stab_\Si(\be,\fb)^\si}
 \mus(\sif)\cn_L(\fb)^{-\fra1{\osf}}\ .}
 Here,
 $[\stab_\Si(\el)\cap\stab_\Si(\be,\fb)\setminus\stab_\Si(\el)]$
 is a set of right coset representatives of
 $\stab_\Si(\el)\cap\stab_\Si(\be,\fb)$ in $\stab_\Si(\el)$. Note
 that the sum $(\al')$ is the sum of all such orbit sums (i).
 \mn Because of \ben\item[]
 $\sif\le\stab_\Si(\be,\fb)^\si\iff\Si_{F^{\si\me}}\le\stab_\Si(\be,\fb)$
 \quad(as $\stab_\Si(\be,\fb)^\si=(\stab_\Si(\be,\fb))^\si$ and $\Si_F^{\si\me}=\Si_{F^{\si\me}}$)\,,
 \item[] $\mus(\sif)=\mus(\Si_{F^{\si\me}})$\quad (a direct consequence of
 the definition of the Möbius function)\,,
 \item[] and $\osf=|\Si_{F^{\si\me}}|$\,,\een
 the inner sums of (i) are independent of $\si$. Hence, if we can
 show
 \nf{ii}{\sum_{F\,\mr{so}\,\sif\le\stab_\Si(\be,\fb)}\mus(\sif)\cn_L(\fb)^{-\fra1{|\sif|}}\equiv0\mod|\stab_\Si(\el)\cap\stab_\Si(\be,\fb)|\zkp\,,}
 then sum (i) is $\equiv0\mod|\stab_\Si(\el)|\zkp$ and the proof of the proposition will be
 complete.
 \mn For (ii), we first shorten the notation by setting $P=\stab_\Si(\be,\fb)\le \Si$
 and $r=\cn_L(\fb)^{-\fra1{|\stab_\Si(\be,\fb)|}}$ which is a unit
 in $\zkp$
 because $\fb=\fa_M\fo_L$ if $\Si_M=P$ (compare with the proof of $(\al')$ above).
  This turns the left hand side of (ii) into
 $\sum_{1\le P'\le P}\mu_P(P')r^{[P:P']}$\,, as obviously
 $\mus(P')=\mu_P(P')$. Applying now the Claim below we obtain
 $$\sum_{1\le P'\le P}\mu_P(P')r^{[P:P']}\equiv 0\mod |P|\zkp\,,$$
 which is even stronger than (ii).
 \mn{\sc Claim\,:} \ {\em Let $P$ be a finite $p$-group and $r$ a unit in $\zkp$. Then}
 $$\sum_{1\le P'\le P}\mu_P(P')r^{[P:P']}\equiv0\mod |P|\zkp\
 .$$
 To see this, let $z\in\zp$
 satisfy $z^{p-1}=1$ and $r\equiv z\mod p$\,, hence $$r^{p^n}\equiv z^{p^n}=z\mod p^{n+1}\,,\quad\mr{for\ }n\ge0\,.$$
 From [HIÖ,p.717] we obtain $|P'|\,\Big|\,p\mu_{P'}(P')$\,, for $P'\le P$\,. Therefore, and as $\mu_P(P')=\mu_{P'}(P')$,
 \ $\mu_P(P')r^{[P:P']}\equiv\mu_P(P')z\mod|P|\zp\, .$
 Consequently,
 $$\sum_{1\le P'\le P}\mu_P(P')r^{[P:P']}\equiv z\cdot\sum_{1\le P'\le
 P}\mu_P(P')=0\mod|P|\zp\ ,$$ as $\sum_{1\le P'\le
 P}\mu_P(P')=0$\,, by the definition of the Möbius function.
 \sn Since $\Q\cap\zp=\zkp$, this proves the Claim
 and also ends the proof of Proposition 9.\qed
 \bn Choosing $\ga\in\hat K\mal$ with $j(\ga)=g$ and denoting by
 $\ee_\ga$ the $q$-expansion of $\ee$ at the cusp determined by $\ga$,
 set $\ee(\ga)\df\cnk(\ga_p)^{-|\Si|k}\ee_\ga$ \ (so $\ee(1)$ is the standard $q$-expansion of $\ee$). Then [DR, (0.3)] implies that
 \bqo{\em the constant term of $\ee(1)-\ee(\ga)$ is contained in
 $p^n\zkp$, provided that $\ee(1)$ has all non-constant
 coefficients contained in $p^n\zkp$\,,}\eqo and, by Proposition 9, this
 applies with \ $p^n=|\stab_\Si(\el)|$\,, so the constant term of
 $\ee(1)-\ee(\ga)$ is in $|\stab_\Si(\el)|\zkp\,.$
 \Section{3}{Conclusion of the proof\,: special cusps}
  By what has been said at the end of the previous
 section we need to compute the constant coefficients of $\ee(\ga)$ and $\ee(1)-\ee(\ga)$\,.
 We do this for the $\ga$'s of Proposition 5 and then prove the Theorem stated in the
 introduction.
 \Lemma{10} Setting $g=j(\ga)$, with $\ga$ as in Proposition 5, the constant term of
 $\ee(\ga)$ is $$\cnk(g)^{|\Si|k}\cdot(\,\sum_F\mus(\sif)\ti\ze_F(1-\osfk,(\el\verfl)_{g_F})\,)\,.$$
 \Stop
 For\,the\,proof\,we\,first\,show\,that $\res_F^KE_F$ has\,constant\,term
 $\cnk((\ga))^{|\Si|k}\ti\ze_F(1-\osfk,\!(\el\verfl)_{g_F})$ at the
 cusp determined by $\ga\in\hat K\mal$.
 By (2) of Lemma 7, this constant term of $\res_F^KE_F$ is equal to the one of $E_F$ at the cusp determined by
 $\ga\in \hat F\mal$\,, whence, by (2) of Proposition 8, equals
 $$\cnf((\ga)_F)^{\osfk}\ti\ze_F(1-\osfk,(\el\ver_F^L)_{g_F})=\cnk((\ga))^{|\Si|k}\ti\ze_F(1-\osfk,(\el\ver_F^L)_{g_F})\,,$$
 because $(\ga)_F=(\ga)_K\fo_F$\,.
 \sn We next check that
 $(\res_F^KE_F)_{|_{|\Si|k}}U_{[F:K]}$ and $\res_F^KE_F$ have the
 same constant term at the cusp determined by $\ga\in \hat K\mal$.
 By {\em a)} of Proposition 5, $M\df\bigl(\begin{smallmatrix}\ga&0\\ 0&\ga\me
 \end{smallmatrix}\bigr)$ is in $\widehat{\gann(\ff)}$\,, so $M=M_1M_2$
 with $M_2=I$ in the notation of [DR,p.262], hence, by [DR,5.8],
 the constant terms referred to above are the constant terms of the standard $q$-expansions of
 $((\res_F^K)_{|_{|\Si|k}}U_{[F:K]})|M_2=(\res_F^K)_{|_{|\Si|k}}U_{[F:K]}$ and
 $(\res_F^KE_F)|M_2=\res_F^KE_F$\,, which agree by Lemma 6.
 \sn Möbius-summing, we have shown that $\ee_\ga$ has constant
 term
 $$\cnk((\ga))^{|\Si|k}\sum_F\mus(\sif)\ti\ze_F(1-|\sif|k,(\el\ver_F^L)_{g_F})\,,$$
 hence $\ee(\ga)$ has the required constant term, because
 $\cnk(\ga_p)\me\cnk((\ga))=\cn_{K,p}(g)$\,, by (3) of Proposition
 8.
 \sn This completes the proof of Lemma 10.\qed
 \mn We can now complete the proof of the Theorem.
 It follows from equation (0) and Lemma 10 that for such $\ga$ the constant term
 of $\ee(1)-\ee(\ga)$ is
 $$\barr{l}\sum_F\mus(\sif)[\ti\ze_F(1-\osfk,\el\verfl)-\cnk(g)^{|\Si|k}\ti\ze_F(1-\osfk,(\el\verfl)_{g_F})]\\[1ex]
 =\sum_F\mus(\sif)\ti\De_{g_F}(1-\osfk,\el\verfl)\  ;\earr$$
 the latter since $\cnf(\verkf g)=\cnk(g)^{[F:K]}$\,.
 \sn Note at this stage that it is this sum, \ $\sum_F\mus(\sif)\ti\De_{g_F}(1-\osfk,\el\verfl)$\,,
 which is referred to in $(4\star)$. The last sentence of \S2 now implies that this sum is
 $\equiv0\mod|\stab_\Si(\el)|\zkp$\,, for every even locally
 constant $\el$, thus verifying the hypothesis of Proposition 4.
 \sn The Theorem is now the conclusion of Proposition 4. \qed
 \remark The application of the Theorem to equivariant Iwasawa theory, in [RW10], is different from that
 proposed in [RW7,\,pp.715/716] since the specialization map $G(L_S/K)\to G({\textbf K}/K)$ extends to a map of
 localizations in the sense of [loc.cit.] only when $G( {\textbf K}/K)$ is a pro-$p$ group.

 \bbn{\sc References}
 \bn\small
 \btb{rp{13cm}}
 \,[CFKSV]& Coates, J., Fukaya, T., Kato, K., Sujatha, R.~and Venjakob, O.\,, {\em The GL$_2$ main
            conjecture for elliptic curves without complex multiplication.} Publ.\,Math.\,Inst.\,Hautes \`Etudes Sci.~{\bf 101}
            (2005), 163-208\\
 \,[DR]&    Deligne, P.~and Ribet, K., {\em Values of abelian
            $L$-functions at negative integers
            over totally real fields.} Invent.~Math. {\bf 59} (1980), 227-286\\
 \,[FK]   & Fukaya, T., Kato, K., {\em A formulation of conjectures on $p$-adic zeta
            functions in noncommutative Iwasawa theory.} Proceedings of
            the St.~Petersburg Mathematical Society, vol.~XII (ed.
            N.N.~Uraltseva), AMS Translations -- Series 2, {\bf 219}
            (2006), 1-86\\
 \,[Ha] &  Hara, T.\,,\newline
           \hsp{2} 1. {\em Iwasawa theory of totally real fields for certain non-commutative $p$-extensions.}\newline
           \hsp{6} arXiv:0803.0211v2 [math.NT]\newline
           \hsp{2} 2. {\em Inductive construction of the $p$-adic zeta functions for non-commutative $p$-\linebreak \hsp{8}extensions
           of totally real fields with exponent $p$.} arXiv:0908.2178v1 [math.NT]\\
 \,[HIÖ]&  Hawkes, T., Isaacs, I.M.~and Özaydin, M.\,, {\em On the Möbius function of a finite group.} Rocky Mountain J.~of
           Mathematics {\bf 19} (1989), 1003-1033 \\
 \,[K] &   Kato, K.\,, {\em Iwasawa theory and generalizations.}
           Proc.\,ICM, Madrid, Spain, 2006; European
           Math.\,Soc.\,(2007), 335-357\\
 \,[Ka] &  Kakde, M.\, {\em Proof of the main conjecture of noncommutative Iwasawa theory for totally real number fields in certain cases.}
           arXiv:0802.2272v2 [math.NT]\\
  \,[RW]&   Ritter, J.~and Weiss, A.\,, \newline
           \hsp{2}2. {\em Towards equivariant Iwasawa theory, II.}
           Indag.\,Mathemat.~{\bf 15} (2004), 549-572\newline
           \hsp{2}3. {\em $\cdots$, III.} Math.\,Ann.~{\bf 336}
           (2006), 27-49\newline
           \hsp{2}7. {\em Congruences between abelian pseudomeasures.}
           Math.Res.\,Lett.~{\bf 15} (2008), 715-725\newline
           \hsp{2}8. {\em Equivariant Iwasawa Theory: An Example.}
           Documenta Math.~{\bf 13} (2008), 117-129\newline
           \hsp{.2}10. {\em On the `main conjecture' of equivariant Iwasawa
           theory.} Preprint (2010\,a)\\
           % arXiv:1004.2578\linebreak \hsp{5} [math.NT]\\
\etb
 \noi
 \btb{rp{13cm}}
 \,[RoW] & Roblot, X.-F.~and Weiss, A.\,, {\em Numerical evidence
           toward a 2-adic equivariant ``main conjecture''.} To appear in
           Experimental Mathematics.\\
 \,[Se]&   Serre, J.-P.\,, {\em Sur le r\'esidu de la fonction
           z\^{e}ta $p$-adique d'un corps de nombres.}\newline
           C.R.Acad.Sci.~Paris {\bf 287} (1978), s\'erie A, 183-188

 \etb
 \tiny
 \vspace*{1.5cm}
 \bct Jürgen Ritter, Schnurbeinstraße 14, 86391 Deuringen,
 Germany; {\tt jr@ritter-maths.de}\\
 Alfred Weiss, Department of Mathematics, University of Alberta,
 Edmonton, AB, Canada T6G 2G1; {\tt weissa@ualberta.ca}\ect

\end{document}